\documentclass[reqno,12pt]{amsart}
\theoremstyle{plain}
\newtheorem{thm}{Theorem}[section]
\theoremstyle{plain}
\newtheorem{cor}[thm]{Corollary} 
\newtheorem{lemma}[thm]{Lemma} 
\newtheorem{prop}[thm]{Proposition}

\newtheorem{example}[thm]{Example}
\newtheorem{examples}[thm]{Examples}

\newtheorem{remark}[thm]{Remark}

\newtheorem{ques}[thm]{Question}

\usepackage{amscd,amssymb,comment,epic,eufrak,euscript,graphics}
\usepackage[all]{xy}

\newcount\theTime
\newcount\theHour
\newcount\theMinute
\newcount\theMinuteTens
\newcount\theScratch
\theTime=\number\time
\theHour=\theTime
\divide\theHour by 60
\theScratch=\theHour
\multiply\theScratch by 60
\theMinute=\theTime
\advance\theMinute by -\theScratch
\theMinuteTens=\theMinute
\divide\theMinuteTens by 10
\theScratch=\theMinuteTens
\multiply\theScratch by 10
\advance\theMinute by -\theScratch

\def\today{{\number\day\space
 \ifcase\month\or
  January\or February\or March\or April\or May\or June\or
  July\or August\or September\or October\or November\or December\fi
 \space\number\year}}

\newcommand\Afr{{\mathfrak A}}

\newcommand\ah{{\hat a}}

\newcommand\alg{{\operatorname{alg}}}

\newcommand\At{{\widetilde A}}

\newcommand\bh{{\hat b}}

\newcommand\Bt{{\widetilde B}}

\newcommand\Cfr{{\mathfrak C}}

\newcommand\Cpx{{\mathbf C}}

\newcommand\Dt{{\widetilde{D}}}

\newcommand\dt{{\tilde{d}}}

\newcommand\HEu{{\EuScript H}}                   

\newcommand\HEut{{\widetilde\HEu}}                   

\newcommand\KEu{{\EuScript K}}                   

\newcommand\phit{{\tilde\phi}}

\newcommand\pit{{\tilde\pi}}

\newcommand\rank{\mathrm{rank}\,}

\newcommand\Reals{{\mathbf R}}

\newcommand\red{{\operatorname{red}}}

\newcommand\xh{{\hat x}}

\begin{document}

\pagestyle{myheadings}

\title{On embeddings of full amalgamated free product C$^*$--algebras \\[1ex] \rm Revised and expanded version}
\author{Scott Armstrong, Ken Dykema, Ruy Exel and Hanfeng Li}

\address{\hskip-\parindent
Scott Armstrong \\
Department of Mathematics \\
University of California \\
Berkeley CA 94720, USA}
\email{sarm@math.berkeley.edu}

\address{\hskip-\parindent
Ken Dykema \\
Department of Mathematics \\
Texas A\&M University \\
College Station TX 77843--3368, USA}
\email{Ken.Dykema@math.tamu.edu}

\address{\hskip-\parindent
Ruy Exel \\
Departamento de Matematica \\
Universidade Federal de Santa Catarina \\
88040-900 Florianopolis SC, BRAZIL}
\email{exel@mtm.ufsc.br}

\address{\hskip-\parindent
Hanfeng Li \\
Department of Mathematics \\
University of Toronto \\
Toronto ON M5S 3G3, CANADA}
\email{hli@fields.toronto.edu}

\thanks{The first author was supported in part by an REU stipend from the NSF.
The second author was supported in part by NSF grant DMS--0070558.
The third author was supported in part by CNPq grant 303968/85--0.
The fourth author was supported jointly by the Mathematics
Department of the University of Toronto and NSERC Grant 8864-02 of
George A. Elliott.}

\date{11 February, 2003}

\begin{abstract}
We examine the question of when the $*$--homomorphism $\lambda:A*_D B\to\At*_\Dt\Bt$
of full amalgamated free product C$^*$--algebras,
arising from compatible inclusions of C$^*$--algebras $A\subseteq\At$, $B\subseteq\Bt$
and $D\subseteq\Dt$, is an embedding.
Results giving sufficient conditions for $\lambda$ to be injective,
as well of classes of examples where $\lambda$ fails
to be injective, are obtained.
As an application, we give necessary and sufficient condition for the full amalgamated free product of
finite dimensional C$^*$--algebras to be residually finite dimensional.
\end{abstract}

\maketitle

\markboth{Full Free Products}{Armstrong, Dykema, Exel, Li}

\section{Introduction}

Given C$^*$--algebras $A$, $B$ and $D$ with injective
$*$--homomorphisms $\phi_A:D\to A$ and $\phi_B:D\to B$, the
corresponding full amalgamated free product C$^*$--algebra
(see~\cite{B78} or~\cite[Chapter\ 5]{Loring}) is the
C$^*$--algebra $\Afr$, equipped with injective $*$--homomorphisms
$\sigma_A:A\to\Afr$ and $\sigma_B:B\to\Afr$ such that
$\sigma_A\circ\phi_A=\sigma_B\circ\phi_B$, such that $\Afr$ is
generated by $\sigma_A(A)\cup\sigma_B(B)$ and satisfying the
universal property that whenever $\Cfr$ is a C$^*$--algebra and
$\pi_A:A\to\Cfr$ and $\pi_B:B\to\Cfr$ are $*$--homomorphisms
satisfying $\pi_A\circ\phi_A=\pi_B\circ\phi_B$, there is a
$*$--homomorphism $\pi:\Afr\to\Cfr$ such that
$\pi\circ\sigma_A=\pi_A$ and $\pi\circ\sigma_B=\pi_B$. This
situation is illustrated by the following commuting diagram:
\begin{equation}\label{eq:defdiag}
\xymatrix{
& D \ar @{_{(}->} [dl]_{\phi_A} \ar @{^{(}->} [dr]^{\phi_B} \\
A\ar @{^{(}->} [r]^{\sigma_A} \ar[dr]_{\pi_A} & {\Afr}\ar @{-->} [d]^--\pi & B.\ar @{_{(}->} [l]_{\sigma_B} \ar[dl]^{\pi_B} \\
& {\Cfr}
}
\end{equation}
The full amalgamated free product C$^*$--algebra $\Afr$ is
commonly denoted by $A*_DB$, although this notation hides the
dependence of $\Afr$ on the embeddings $\phi_A$ and $\phi_B$.

\begin{ques}\label{ques:emb}
Let $D$, $A$, $B$, $\Dt$, $\At$ and $\Bt$ be C$^*$--algebras and suppose there are injective $*$--homomorphisms
making the following diagram commute:
\[
\xymatrix{
{\At} & {\Dt} \ar @{_{(}->} [l]_{\phi_\At} \ar @{^{(}->} [r]^{\phi_\Bt} & {\Bt.} \\
A \ar @{^{(}->} [u]^{\lambda_A} & D \ar @{_{(}->} [l]_{\phi_A} \ar @{^{(}->} [u]^{\lambda_D} \ar @{^{(}->} [r]^{\phi_B}
 & B \ar @{^{(}->} [u]^{\lambda_B}
}
\]
Let $A*_DB$ and $\At*_\Dt\Bt$ be the corresponding full amalgamated free product C$^*$--algebras
and let $\lambda:A*_DB\to\At*_\Dt\Bt$ be the $*$--homomorphism arising from $\lambda_A$ and $\lambda_B$ via the universal
property.
When is $\lambda$ injective?
\end{ques}

We prove in~\S\ref{sec:emb} that $\lambda$ is injective
when either (i)~$D=\Dt$, (or more precisely, when the $*$--homomorphism $\lambda_D$
is surjective), or (ii)~there are conditional expectations $E_A:\At\to A$ and $E_B:\Bt\to B$ that send $\Dt$ onto $D$ and agree on $\Dt$.
Injectivity in the case $D=\Dt$ was previously proved by G.K.\ Pedersen~\cite{Pedersen}.
(Moreover, earlier results of F.~Boca~\cite{Boca91}
imply that the map $\lambda$ is injective when $D=\Dt$ and when there are conditional
expectations
\[
\At\overset{E^\At_A}{\rightarrow}A\overset{E^A_D}{\rightarrow}D\overset{E^B_D}{\leftarrow}B\overset{E^\Bt_B}{\leftarrow}\Bt;
\]
an argument for the case $D=\Dt=\Cpx$, which uses Boca's results, is outlined in~\cite[4.7]{BP}.)
However, we include our proof because it is different from that found in~\cite{Pedersen}
and because it contains the main idea of our proof of injectivity in case~(ii).
In~\S\ref{sec:ex}, we consider some general conditions and give some concrete examples when $\lambda$ fails to be injective.
Finally, in~\S\ref{sec:rfd}, we apply this embedding result to extend a result from~\cite{BD} about residual finite dimensionality
of full amalgamated free products of finite dimensional C$^*$--algebras.

\section{Embeddings of full free products}\label{sec:emb}

The following result is of course well known.
We include a proof for completeness.
\begin{lemma}\label{lem:repextend}
Let $A$ be a C$^*$--subalgebra of a C$^*$--algebra $\At$ and let $\pi:A\to B(\HEu)$ be a $*$--representation.
Then there is a Hilbert space $\KEu$ and a $*$--representation $\pit:\At\to B(\HEu\oplus\KEu)$ such that
\begin{equation}\label{eq:repextend}
\pit(a)(h\oplus0)=(\pi(a)h)\oplus0,\qquad(a\in A,\,h\in\HEu).
\end{equation}
\end{lemma}
\begin{proof}
Since in general $\pi$ is a direct sum of cyclic representations, we may without loss of generality assume $\pi$
is a cyclic representation with cyclic vector $\xi$.
Let $\phi$ be the vector state $\phi(\cdot)=\langle\pi(\cdot)\xi,\xi\rangle$ of $A$.
Then $\HEu$ is identified with $L^2(A,\phi)$ and $\pi$ is the associated GNS representation.
Let $\phit$ be an extension of $\phi$ to a state of $\At$ and let $\HEut=L^2(\At,\phit)$.
Then the inclusion $A\hookrightarrow\At$ gives rise to an isometry $\HEu\to\HEut$, and we may thus write $\HEut=\HEu\oplus\KEu$
for a Hilbert space $\KEu$.
If $\pit:\At\to B(\HEu\oplus\KEu)$ is the GNS representation associated to $\phit$, then~\eqref{eq:repextend} holds.
\end{proof}

The following result was first proved by G.K.\ Pedersen~\cite{Pedersen}.
We offer a new proof, which is perhaps more elementary.
This proof contains essentially the
same idea as our proof of Proposition~\ref{prop:condexpembed} below.

\begin{prop}\label{prop:fpembed}
Let
\[
\At\supseteq A\supseteq D\subseteq B\subseteq\Bt
\]
be inclusions of C$^*$--algebras and let $A*_DB$ and $\At*_D\Bt$ be
the corresponding full amalgamated free product C$^*$--algebras.
Let $\lambda:A*_DB\to\At*_D\Bt$ be the $*$--homomorphism arising via the universal property from the inclusions $A\hookrightarrow\At$ and
$B\hookrightarrow\Bt$.
Then $\lambda$ is injective.
\end{prop}
\begin{proof}
Let $\pi:A*_DB\to B(\HEu)$ be a faithful $*$--homomorphism.
We will find a Hilbert space $\KEu$ and a $*$--homomorphism $\pit:\At*_D\Bt\to B(\HEu\oplus\KEu)$ such that
\begin{equation}\label{eq:pit}
\pit(\lambda(x))(h\oplus0)=(\pi(x)h)\oplus0,\qquad(x\in A*_DB,\,h\in\HEu).
\end{equation}
This will imply $\lambda$ is injective.

Let $\pi_A:A\to B(\HEu)$ and $\pi_B:B\to B(\HEu)$ be the $*$--representations obtained by composing $\pi$ with the
inclusions $A\hookrightarrow A*_DB$ and $B\hookrightarrow A*_DB$.
Let
\begin{align*}
\sigma_{A,0}&:\At\to B(\HEu\oplus\KEu_{A,0}), \\
\sigma_{B,0}&:\Bt\to B(\HEu\oplus\KEu_{B,0})
\end{align*}
be $*$--representations obtained from Lemma~\ref{lem:repextend} such that
\begin{equation*}
\sigma_{A,0}(a)(h\oplus0)=(\pi_A(a)h)\oplus0,\qquad(a\in A,\,h\in\HEu),
\end{equation*}
and similarly with $A$ replaced by $B$.
Note that $0\oplus\KEu_{A,0}$ is reducing for $\sigma_{A,0}(D)$.
Using Lemma~\ref{lem:repextend},
find Hilbert spaces $\KEu_{B,1}$ and $\KEu_{A,1}$ and $*$--representations
\begin{align*}
\sigma_{B,1}:&\Bt\to B(\KEu_{A,0}\oplus\KEu_{B,1}) \\
\sigma_{A,1}:&\At\to B(\KEu_{B,0}\oplus\KEu_{A,1})
\end{align*}
such that
\begin{align*}
\sigma_{B,1}(d)(k\oplus0)=\sigma_{A,0}(d)(0\oplus k),\qquad(d\in D,\,k\in\KEu_{A,0}), \\
\sigma_{A,1}(d)(k\oplus0)=\sigma_{B,0}(d)(0\oplus k),\qquad(d\in D,\,k\in\KEu_{B,0}).
\end{align*}
Proceeding recursively, for every integer $n\ge2$ we find $*$--representations
\begin{align*}
\sigma_{B,n}&:\Bt\to B(\KEu_{A,n-1}\oplus\KEu_{B,n}), \\
\sigma_{A,n}&:\Bt\to B(\KEu_{B,n-1}\oplus\KEu_{A,n})
\end{align*}
such that
\begin{align*}
\sigma_{B,n}(d)(k\oplus0)&=\sigma_{A,n-1}(d)(0\oplus k),\qquad(d\in D,\,k\in\KEu_{A,n-1}), \\
\sigma_{A,n}(d)(k\oplus0)&=\sigma_{B,n-1}(d)(0\oplus k),\qquad(d\in D,\,k\in\KEu_{B,n-1}).
\end{align*}

We now define the Hilbert spaces
\begin{equation}\label{eq:HEut}
\begin{aligned}
\HEut_A&=\overset{\sigma_{A,0}}{\overbrace{\HEu\oplus\KEu_{A,0}}}
 \oplus\overset{\sigma_{A,1}}{\overbrace{\KEu_{B,0}\oplus\KEu_{A,1}}}
 \oplus\overset{\sigma_{A,2}}{\overbrace{\KEu_{B,1}\oplus\KEu_{A,2}}}\oplus\cdots, \\
\HEut_B&=\underset{\sigma_{B,0}}{\underbrace{\HEu\oplus\KEu_{B,0}}}
 \oplus\underset{\sigma_{B,1}}{\underbrace{\KEu_{A,0}\oplus\KEu_{B,1}}}
 \oplus\underset{\sigma_{B,2}}{\underbrace{\KEu_{A,1}\oplus\KEu_{B,2}}}\oplus\cdots,
\end{aligned}
\end{equation}
where the brackets indicate where the constructed representations act,
and we let $\sigma_\At:\At\to B(\HEut_A)$ and $\sigma_\Bt:\Bt\to B(\HEut_B)$ be the $*$--representations
\begin{align*}
\sigma_\At&=\sigma_{A,0}\oplus\sigma_{A,1}\oplus\sigma_{A,2}\oplus\cdots, \\
\sigma_\Bt&=\sigma_{B,0}\oplus\sigma_{B,1}\oplus\sigma_{B,2}\oplus\cdots,
\end{align*}
where the summands act as indicated by brackets in~\eqref{eq:HEut}.
Consider the unitary $U:\HEut_A\to\HEut_B$ mapping the summands in $\HEut_A$ identically to the corresponding
summands in $\HEut_B$ as indicated by the arrows below:
\[
\xymatrix@C=0em{
\HEut_A\ar[d]_{U} & = & \HEu \ar[d] & \oplus & \KEu_{A,0} \ar[drr] & \oplus & \KEu_{B,0} \ar[dll] & \oplus
  & \KEu_{A,1} \ar[drr] & \oplus & \KEu_{B,1} \ar[dll] & \oplus & \KEu_{A,2} \ar[drr] & \oplus & \cdots \ar[dll] \\
\HEut_B    & = & \HEu   & \oplus & \KEu_{B,0} & \oplus & \KEu_{A,0} & \oplus
  & \KEu_{B,1} & \oplus & \KEu_{A,1} & \oplus & \KEu_{B,2} & \oplus & \cdots \;.
}
\]
Let $\KEu=\KEu_{A,0}\oplus\KEu_{B,0}\oplus\KEu_{A,1}\oplus\KEu_{B,1}\oplus\cdots$ and identify $\HEu\oplus\KEu$ with $\HEut_A$.
Then we have the $*$--representations $\pit_\At=\sigma_\At:\At\to B(\HEu\oplus\KEu)$ and
$\pit_\Bt:\Bt\to B(\HEu\oplus\KEu)$, the latter defined by $\pit_\Bt(\cdot)=U^*\sigma_\Bt(\cdot)U$.
By construction, the restrictions of $\pit_\At$ and $\pit_\Bt$ to $D$ agree, and we have
\begin{align*}
\pit_\At(a)(h\oplus0)&=(\pi_A(a)h)\oplus0,\qquad(a\in A,\,h\in\HEu), \\
\pit_\Bt(b)(h\oplus0)&=(\pi_B(b)h)\oplus0,\qquad(b\in B,\,h\in\HEu).
\end{align*}
Letting $\pit:\At*_D\Bt\to B(\HEu\oplus\KEu)$ be the $*$--homomorphism obtained from $\pit_\At$ and $\pit_\Bt$ via the
universal property, we have that~\eqref{eq:pit} holds.
\end{proof}

For a C$^*$--algebra $A$, unital or not, let $A^u$ denote the unitization of $A$.
Thus, as a vector space, $A^u=A\oplus\Cpx$ with multiplication defined by
$(a,\mu)\cdot(a',\mu')=(aa'+\mu a+\mu' a,\mu\mu')$.
We identify $A$ with the ideal $A\oplus0$ of $A^u$, which has codimension $1$.

\begin{lemma}\label{lem:unitization}
Let $A\supseteq D\subseteq B$ be inclusions of C$^*$--algebras.
Consider the unitizations and corresponding inclusions
\[
\xymatrix{
{A^u} & {D^u} \ar @{_{(}->} [l] \ar @{^{(}->} [r] & {B^u} \\
A \ar @{^{(}->} [u] & D \ar @{_{(}->} [l] \ar @{^{(}->} [u] \ar @{^{(}->} [r]
 & {B.} \ar @{^{(}->} [u]
}
\]
Let $\lambda:A*_DB\to A^u*_{D^u}B^u$ be the resulting $*$--homomorphism between full amalgamated free products.
Then there is an isomorphism $\pi:A^u*_{D^u}B^u\to(A*_DB)^u$ such that $\pi\circ\lambda:A*_DB\to(A*_DB)^u$ is the 
canonical embedding arising in the definition of the unitization.
\end{lemma}
\begin{proof}
Since any $*$--representations of $A$ and $B$ that agree on $D$ extend to $*$--rep\-resen\-tations of $A^u$ and $B^u$
that agree on $D^u$, the $*$--homomorphism $\lambda$ is injective.
Let $e\in A^u*_{D^u}B^u$ be the unit of $A^u$, which is of course identified with the units of $B^u$ and $D^u$.
Clearly, $A^u*_{D^u}B^u$ is generated by the image of $\lambda$ together with $e$.
One easily sees
\[
(\lambda(x)+\mu e)(\lambda(x')+\mu'e)=\lambda(xx')+\mu\lambda(x')+\mu'\lambda(x)+\mu\mu'e.
\]
Moreover, if $\rho:A^u*_{D^u}B^u\to\Cpx$ is the $*$--homomorphism arising from the unital $*$--homomorphisms $A^u\to\Cpx$
and $B^u\to\Cpx$, then $\rho(e)=1$ and $\lambda(A*_DB)\subseteq\ker\rho$.
Hence $\lambda(A*_D B)$ has codimension $1$ in $A^u*_{D^u}B^u$.
Now $\pi$ can be defined by $\pi(\lambda(x)+\mu e)=(x,\mu)$.
\end{proof}

\begin{prop}\label{prop:condexpembed}
Suppose
\begin{equation}\label{eq:incl}
\xymatrix{
{\At} & {\Dt} \ar @{_{(}->} [l] \ar @{^{(}->} [r] & {\Bt} \\
A \ar @{^{(}->} [u] & D \ar @{_{(}->} [l] \ar @{^{(}->} [u] \ar @{^{(}->} [r]
 & B \ar @{^{(}->} [u]
}
\end{equation}
is a commuting diagram of inclusions of C$^*$--algebras.
Let $\lambda:A*_D B\to\At*_\Dt\Bt$ be the resulting $*$--homomorphism of full free product C$^*$--algebras.
Suppose there are conditional
expectations $E_A:\At\to A$, $E_D:\Dt\to D$ and $E_B:\Bt\to B$ onto $A$, $D$ and $B$,
respectively, such that the diagram
\begin{equation}\label{eq:condexp}
\xymatrix{
{\At} \ar [d]^{E_A} & {\Dt} \ar [d]^{E_D} \ar @{_{(}->} [l] \ar @{^{(}->} [r] & {\Bt} \ar [d]^{E_B} \\
A & D \ar @{_{(}->} [l] \ar @{^{(}->} [r] & B
}
\end{equation}
commutes.
Then $\lambda$ is injective.
\end{prop}
\begin{proof}
By appealing to Lemma~\ref{lem:unitization}, we may without loss of generality assume
all the algebras and $*$--homomorphisms in~\eqref{eq:incl} are unital.
Let $\pi:A*_DB\to B(\HEu)$ be a faithful, unital $*$--representation.
As in the proof of Proposition~\ref{prop:fpembed}, in order to show $\lambda$ is injective,
we will find a Hilbert space $\KEu$ and a $*$--homomorphism $\pit:\At*_\Dt\Bt\to B(\HEu\oplus\KEu)$ such that
\begin{equation}\label{eq:pit2}
\pit(\lambda(x))(h\oplus0)=(\pi(x)h)\oplus0,\qquad(x\in A*_DB,\,h\in\HEu).
\end{equation}

Let $\pi_A:A\to B(\HEu)$ and $\pi_B:B\to B(\HEu)$ be the $*$--representations obtained by composing $\pi$ with the
inclusions $A\hookrightarrow A*_DB$ and $B\hookrightarrow A*_DB$,
and let $\pi_D:D\to B(\HEu)$ be their common restriction to $D$.
Consider the
canonical left action of $\Dt$ on the right Hilbert $D$--module $L^2(\Dt,E_D)$, which is obtained from $\Dt$ by separation 
and completion with respect to the $D$--valued inner product $\langle\dt_1,\dt_2\rangle=E_D(\dt_1^*\dt_2)$.
Consider the Hilbert space $L^2(\Dt,E_D)\otimes_D\HEu$, where the left action of $D$ on $\HEu$ is via $\pi_D$.
Since $\pi_D$ is unital, $\HEu$ embeds as a subspace, and we can write 
\begin{equation}\label{eq:L2D}
L^2(\Dt,E_D)\otimes_D\HEu=\HEu\oplus\KEu_D.
\end{equation}
Consider the left action of $\Dt$
on the Hilbert space $\HEu\oplus\KEu_D$.
The subspace $\HEu$ is reducing for the restriction of $\sigma_D$ to $D$, and we have $\sigma_D(d)(h\oplus0)=(\pi_D(d)h)\oplus0$ for
every $d\in D$ and $h\in\HEu$.

In a similar way, consider the Hilbert spaces
\begin{equation}\label{eq:L2AB}
L^2(\At,E_A)\otimes_A\HEu,\qquad L^2(\Bt,E_B)\otimes_B\HEu
\end{equation}
and the associated left actions $\sigma_{A,0}$ of $\At$, respectively $\sigma_{B,0}$ of $\Bt$.
As the diagram~\eqref{eq:condexp} commutes, the Hilbert space~\eqref{eq:L2D} embeds canonically
as a subspace of both spaces~\eqref{eq:L2AB}.
We may thus write
\begin{align*}
L^2(\At,E_A)\otimes_A\HEu&=\HEu\oplus\KEu_D\oplus\KEu_{A,0} \\
L^2(\Bt,E_B)\otimes_B\HEu&=\HEu\oplus\KEu_D\oplus\KEu_{B,0},
\end{align*}
the subspace $\HEu\oplus\KEu_D\oplus0$ is reducing for the restrictions
of $\sigma_{A,0}$ and $\sigma_{B,0}$ to $\Dt$, and we have
$\sigma_{A,0}(\dt)(\eta\oplus0)=(\sigma_D(\dt)\eta)\oplus0=\sigma_{B,0}(\dt)(\eta\oplus0)$
for every $\dt\in\Dt$ and $\eta\in\HEu\oplus\KEu_D$.
Moreover, $\HEu\oplus0\oplus0$ is reducing for the restrictions of $\sigma_{A,0}$ to $A$
and $\sigma_{B,0}$ to $B$, and we have
\begin{align*}
\sigma_{A,0}(a)(h\oplus0\oplus0)&=(\pi_A(a)h)\oplus0\oplus0\qquad(a\in A,\,h\in\HEu) \\
\sigma_{B,0}(b)(h\oplus0\oplus0)&=(\pi_B(b)h)\oplus0\oplus0\qquad(b\in B,\,h\in\HEu).
\end{align*}
Let $\sigma_{A,0,\Dt}$ denote the action of $\Dt$ on $\KEu_{A,0}$ obtained by restricting $\sigma_{A,0}$ to $\Dt$ and compressing,
and similarly for $\sigma_{B,0,\Dt}$.

We now proceed recursively as in the proof of Proposition~\ref{prop:fpembed}.
If Hilbert spaces $\KEu_{A,n-1}$ and $\KEu_{B,n-1}$ have been constructed with actions $\sigma_{A,n-1,\Dt}$
and $\sigma_{B,n-1,\Dt}$, respectively, of $\Dt$,
use Lemma~\ref{lem:repextend} to construct Hilbert spaces $\KEu_{B,n}$ and $\KEu_{A,n}$ and $*$--homomorphisms
\begin{align*}
\sigma_{B,n}:\Bt&\to B(\KEu_{A,n-1}\oplus\KEu_{B,n}) \\
\sigma_{A,n}:\At&\to B(\KEu_{B,n-1}\oplus\KEu_{A,n}),
\end{align*}
such that
\begin{align*}
\sigma_{B,n}(\dt)(k\oplus0)&=(\sigma_{A,n-1,\Dt}(\dt)k)\oplus0\qquad(\dt\in\Dt,\,k\in\KEu_{A,n-1}) \\
\sigma_{A,n}(\dt)(k\oplus0)&=(\sigma_{B,n-1,\Dt}(\dt)k)\oplus0\qquad(\dt\in\Dt,\,k\in\KEu_{B,n-1}).
\end{align*}
Then let $\sigma_{B,n,\Dt}$ be the action of $\Dt$ on $\KEu_{B,n}$ obtained from the restriction of $\sigma_{B,n}$
to $\Dt$ by compressing, and similarly define the action $\sigma_{A,n,\Dt}$ of $\Dt$ on $\KEu_{A,n}$.

We may now define the Hilbert spaces
\begin{equation}\label{eq:HEut2}
\begin{aligned}
\HEut_A&=\overset{\sigma_{A,0}}{\overbrace{\overset{\sigma_D}{\overbrace{\HEu\oplus\KEu_D}}\oplus\KEu_{A,0}}}
 \oplus\overset{\sigma_{A,1}}{\overbrace{\KEu_{B,0}\oplus\KEu_{A,1}}}
 \oplus\overset{\sigma_{A,2}}{\overbrace{\KEu_{B,1}\oplus\KEu_{A,2}}}\oplus\cdots, \\[2ex]
\HEut_B&=\underset{\sigma_{B,0}}{\underbrace{\underset{\sigma_D}{\underbrace{\HEu\oplus\KEu_D}}\oplus\KEu_{B,0}}}
 \oplus\underset{\sigma_{B,1}}{\underbrace{\KEu_{A,0}\oplus\KEu_{B,1}}}
 \oplus\underset{\sigma_{B,2}}{\underbrace{\KEu_{A,1}\oplus\KEu_{B,2}}}\oplus\cdots,
\end{aligned}
\end{equation}
where the brackets indicate where the constructed representations act.
We let $\sigma_\At:\At\to B(\HEut_A)$ and $\sigma_\Bt:\Bt\to B(\HEut_B)$ be the $*$--representations
\begin{align*}
\sigma_\At&=\sigma_{A,0}\oplus\sigma_{A,1}\oplus\sigma_{A,2}\oplus\cdots, \\
\sigma_\Bt&=\sigma_{B,0}\oplus\sigma_{B,1}\oplus\sigma_{B,2}\oplus\cdots,
\end{align*}
where the summands act as indicated by brackets in~\eqref{eq:HEut2}.
Consider the unitary $U:\HEut_A\to\HEut_B$ mapping the summands in $\HEut_A$ identically to the corresponding
summands in $\HEut_B$ as indicated by the arrows below:
\[
\xymatrix@C=0em{
\HEut_A\ar[d]_{U} & = & \HEu \ar[d] & \oplus & \KEu_D \ar[d] & \oplus & \KEu_{A,0} \ar[drr] & \oplus & \KEu_{B,0} \ar[dll] & \oplus
  & \KEu_{A,1} \ar[drr] & \oplus & \KEu_{B,1} \ar[dll] & \oplus & \KEu_{A,2} \ar[drr] & \oplus & \cdots \ar[dll] \\
\HEut_B    & = & \HEu & \oplus & \KEu_D  & \oplus & \KEu_{B,0} & \oplus & \KEu_{A,0} & \oplus
  & \KEu_{B,1} & \oplus & \KEu_{A,1} & \oplus & \KEu_{B,2} & \oplus & \cdots \;.
}
\]
Let $\KEu=\KEu_D\oplus\KEu_{A,0}\oplus\KEu_{B,0}\oplus\KEu_{A,1}\oplus\KEu_{B,1}\oplus\cdots$ and identify $\HEu\oplus\KEu$ with $\HEut_A$.
Then we have the $*$--representations $\pit_\At=\sigma_\At:\At\to B(\HEu\oplus\KEu)$ and
$\pit_\Bt:\Bt\to B(\HEu\oplus\KEu)$, the latter defined by $\pit_\Bt(\cdot)=U^*\sigma_\Bt(\cdot)U$.
By construction, the restrictions of $\pit_\At$ and $\pit_\Bt$ to $\Dt$ agree, and we have
\begin{align*}
\pit_\At(a)(h\oplus0)&=(\pi_A(a)h)\oplus0\qquad(a\in A,\,h\in\HEu) \\
\pit_\Bt(b)(h\oplus0)&=(\pi_B(b)h)\oplus0\qquad(b\in B,\,h\in\HEu).
\end{align*}
Letting $\pit:\At*_\Dt\Bt\to B(\HEu\oplus\KEu)$ be the $*$--homomorphism obtained from $\pit_\At$ and $\pit_\Bt$ via the
universal property, we have that~\eqref{eq:pit2} holds.
\end{proof}

\section{Examples of non--embedding}\label{sec:ex}

In this section, we give some examples when the map $\lambda$ of Question~\ref{ques:emb} fails to be injective.
(In contrast, it is known~\cite{BlDy} that in the more stringent situation of
{\em reduced} amalgamated free products, the map analogous to $\lambda$ is always injective.)

We begin with a trivial class of examples.
\begin{examples}\label{exs:triv}\rm
Let $A$ and $B$ be C$^*$-subalgebras of a
C$^*$-algebra $E$ with $A \nsubseteq B$ and $B\nsubseteq A$. Let
$D=A\cap B, \, \At=E$ and $\Dt=\Bt=B$, equipped with the natural
inclusions. Then the map $\lambda: A*_DB\rightarrow
\At*_{\Dt}\Bt=E$ is injective if and only if $A*_DB$ is exactly
the C$^*$-subalgebra of $E$ generated by $A$ and $B$. This doesn't
hold in general. Notice that in these examples,
$B\cap\Dt=B\supsetneqq D$.
\end{examples}

\begin{prop}\label{prop:ADBnoninj}
Suppose
\begin{equation*}
\xymatrix{
{\At} & {\Dt} \ar @{_{(}->} [l] \ar @{^{(}->} [r] & {\Bt} \\
A \ar @{^{(}->} [u] & D \ar @{_{(}->} [l] \ar @{^{(}->} [u] \ar @{^{(}->} [r]
 & B \ar @{^{(}->} [u]
}
\end{equation*}
is a commuting diagram of inclusions of C$^*$--algebras
and let $\lambda:A*_D B\to\At*_\Dt\Bt$ be the resulting $*$--homomorphism of full free product C$^*$--algebras.
Suppose there are conditional
expectations $E^A_D:A\to D$ and $E^B_D:B\to D$ with $E^B_D$ faithful.
Suppose there are $\dt\in\Dt$, $a\in A$ and $b\in B$ satisfying
$a\dt\in A$, $\dt b\in B$,
\begin{align}
D(\dt b)\cap Db&=\{0\} \label{eq:Db} \\
E^A_D(\dt^*a^*ad)b&\ne0. \label{eq:Eda}
\end{align}
Then $\lambda$ is not injective.
\end{prop}
\begin{proof}
Letting
\begin{equation}\label{eq:sigmas}
\begin{alignedat}{2}
\sigma_A&:A\hookrightarrow A*_DB,\qquad&\sigma_B&:B\hookrightarrow A*_DB, \\
\sigma_\At&:\At\hookrightarrow\At*_\Dt\Bt,&\sigma_\Bt&:\Bt\hookrightarrow\At*_\Dt\Bt
\end{alignedat}
\end{equation}
be the embeddings as in~\eqref{eq:defdiag}, we have
\[
\lambda(\sigma_A(a\dt)\sigma_B(b))=\sigma_\At(a\dt)\sigma_\Bt(b)=\sigma_\At(a)\sigma_\Bt(\dt b)=\lambda(\sigma_A(a)\sigma_B(\dt b)).
\]
Thus we need only show
\begin{equation}\label{eq:sigmane}
\sigma_A(a\dt)\sigma_B(b)\ne\sigma_A(a)\sigma_B(\dt b).
\end{equation}
We consider the reduced amalgamated free product of C$^*$--algebras (see~\cite{Vo} or~\cite{VDN}),
\[
(A*_D^\red B,E_D)=(A,E^A_D)*_D(B,E^B_D)
\]
and the natural quotient $*$--homomorphism $A*_DB\to A*_D^\red B$.
Let $L^2(A*_D^\red B,E_D)$ be the right Hilbert $D$--module obtained by separation and completion from
$A*_D^\red B$ with respect to the $D$--valued inner product $\langle x,y\rangle=E_D(x^*y)$,
and given $x\in A*_D^\red B$, let $\xh$ denote the corresponding element in $L^2(A*_D^\red B,E_D)$.
Let $\HEu_A=L^2(A,E^A_D)$ and $\HEu_B=L^2(B,E^B_D)$ be similarly defined.
Then in $L^2(A*_D^\red B,E_D)$, the closure of the subspace spanned by elements of the form $(ab)\hat{\;}$
for $a\in A$ and $b\in B$ is isomorphic to the tensor product $\HEu_A\otimes_D\HEu_B$ of Hilbert $D$--modules.
In order to show~\eqref{eq:sigmane}, it will suffice to show
\[
(a\dt)\hat{\;}\otimes\bh\ne\ah\otimes(\dt b)\hat{\;}
\]
in $\HEu_A\otimes_D\HEu_B$.
Let $\zeta_B\in\HEu_B$.
Then
\begin{align}
\langle(a\dt)\hat{\;}\otimes\zeta_B,(a\dt)\hat{\;}\otimes\bh\rangle=\langle\zeta_B,(E_D^A(\dt^*a^*a\dt)b)\hat{\;}\rangle
\label{eq:adtb} \\
\langle(a\dt)\hat{\;}\otimes\zeta_B,\ah\otimes(\dt b)\hat{\;}\rangle=\langle\zeta_B,(E_D^A(\dt^*a^*a)\dt b)\hat{\;}\rangle.
\label{eq:atdb}
\end{align}
From assumptions~\eqref{eq:Db} and~\eqref{eq:Eda}, we obtain $E_D^A(\dt^*a^*a\dt)b\ne E_D^A(\dt^*a^*a)\dt b$.
Since $E^B_D$ is faithful, there is $\zeta_B\in\HEu_B$ such that the right--hand--sides of~\eqref{eq:adtb} and~\eqref{eq:atdb}
are not equal.
\end{proof}

\begin{remark}\label{rmk:wkhyp}\rm
From the above proof, one sees that the hypotheses of Proposition~\ref{prop:ADBnoninj}
can be weakened as follows:
Assumptions~\eqref{eq:Db} and~\eqref{eq:Eda} can be dropped, and $E^B_D$ need not be assumed faithful, but instead
one must assume
\begin{equation}\label{eq:Econd}
E^B_D\big(b^*\big(E^A_D(\dt^*a^*a\dt)-E^A_D(\dt^*a^*a)\dt-\dt^*E^A_D(a^*a\dt)+\dt^*E^A_D(a^*a)\dt\big)b^*\big)\ne0.
\end{equation}
Note that the LHS of~\eqref{eq:Econd} is nothing other than
\[
\langle(a\dt)\hat{\;}\otimes\bh-\ah\otimes(\dt b)\hat{\;},(a\dt)\hat{\;}\otimes\bh-\ah\otimes(\dt b)\hat{\;}\rangle.
\]
\end{remark}

\begin{cor}\label{cor:ABnoninj}
Suppose
\begin{equation}\label{eq:incl2}
\xymatrix{
{\At} & {\Dt} \ar @{_{(}->} [l] \ar @{^{(}->} [r] & {\Bt} \\
A \ar @{^{(}->} [u] & D \ar @{_{(}->} [l] \ar @{^{(}->} [u] \ar @{^{(}->} [r]
 & B \ar @{^{(}->} [u]
}
\end{equation}
is a commuting diagram of inclusions of C$^*$--algebras
and let $\lambda:A*_D B\to\At*_\Dt\Bt$ be the resulting $*$--homomorphism of full free product C$^*$--algebras.
Suppose one of the following holds:
\renewcommand{\labelenumi}{(\roman{enumi})}
\begin{enumerate}
\item $D=0$
\item $D=\Cpx$, $A$ and $B$ are unital and the inclusions $D\hookrightarrow A$ and $D\hookrightarrow B$ are unital.
\end{enumerate}
Suppose there are $\dt\in\Dt$, $a\in A$ and $b\in B$ such that $a\dt\in A\backslash\{0\}$, $\dt b\in B$ and $\dt b\notin\Cpx b$.
Then $\lambda$ is not injective.
\end{cor}
\begin{proof}
We can reduce to the case in which~(ii) holds by application of Lemma~\ref{lem:unitization}.
We may without loss of generality assume $A$ and $B$ are separable.
Letting $E^A_D:A\to\Cpx$ and $E^B_D:B\to\Cpx$ be faithful states, we find the hypotheses of Proposition~\ref{prop:ADBnoninj}
are satisfied.
\end{proof}

From this corollary, we have the following class of concrete examples, which
shows that $\lambda$ may be non--injective even if
\begin{equation}\label{eq:ABDtneD}
B\cap\Dt=D=A\cap\Dt.
\end{equation}

\begin{example}\label{ex:noninjective}\rm
Let $\HEu$ be an infinite dimensional, separable Hilbert space.
Inside $B(\HEu)$, let $D=\Cpx1$ and let $A=B=D+K(\HEu)$, where $K(\HEu)$ is the compact operators.
Let $u\in B(\HEu)$ be a unitary operator that does not belong to $D$
and let $\Dt=C^*(u)$, $\At=\Bt=\Dt+K(\HEu)$.
Let $\lambda:A*_DB\to\At*_\Dt\Bt$ be the $*$--homomorphism arising from the inclusions~\eqref{eq:incl2}.
Then $\lambda$ is not injective.
\end{example}
\begin{proof}
Take $\dt=u$ and $a\in K(\HEu)\backslash\{0\}$.
Since $u\notin\Cpx1$, there is $b\in K(\HEu)$ such that $ub\notin\Cpx b$.
Now apply Corollary~\ref{cor:ABnoninj}.
One can choose $u$ so that $C^*(u)\cap (\Cpx1+K(\HEu))=\Cpx1$, in order to get~\eqref{eq:ABDtneD}.
\end{proof}

\begin{prop}\label{prop:db=bd}
Suppose
\begin{equation}\label{eq:incl3}
\xymatrix{
{\At} & {\Dt} \ar @{_{(}->} [l] \ar @{^{(}->} [r] & {\Bt} \\
A \ar @{^{(}->} [u] & D \ar @{_{(}->} [l] \ar @{^{(}->} [u] \ar @{^{(}->} [r]
 & B \ar @{^{(}->} [u]
}
\end{equation}
is a commuting diagram of inclusions of C$^*$--algebras
and let $\lambda:A*_D B\to\At*_\Dt\Bt$ be the resulting $*$--homomorphism of full free product C$^*$--algebras.
Suppose one of the following holds:
\renewcommand{\labelenumi}{(\roman{enumi})}
\begin{enumerate}
\item $D=0$
\item $D=\Cpx$, $A$ and $B$ are unital and the inclusions $D\hookrightarrow A$ and $D\hookrightarrow B$ are unital.
\end{enumerate}
Suppose there are $\dt\in\Dt$, $a_1,a_2\in A$ and $b\in B\backslash D$ such that $a_1\dt,\dt a_2\in A$, $a_1\dt\notin\Cpx$
and $\dt b=b\dt$.
Then $\lambda$ is not injective.
\end{prop}
\begin{proof}
We can reduce to the case in which~(ii) holds by application of Lemma~\ref{lem:unitization}.
We use the same notation as in~\eqref{eq:sigmas}.
We have
\begin{align*}
\lambda(\sigma_A(a_1\dt)\sigma_B(b)\sigma_A(a_2))&=\sigma_\At(a_1\dt)\sigma_\Bt(b)\sigma_\At(a_2) \\
&=\sigma_\At(a_1)\sigma_\Bt(b)\sigma_\At(\dt a_2)=\lambda(\sigma_A(a_1)\sigma_B(b)\sigma_A(\dt a_2)),
\end{align*}
and we must only show
\begin{equation}\label{eq:adbne}
\sigma_A(a_1\dt)\sigma_B(b)\sigma_A(a_2)\ne\sigma_A(a_1)\sigma_B(b)\sigma_A(\dt a_2).
\end{equation}
Without loss of generality, assume $A$ and $B$ are separable.
Let $\phi_A:A\to\Cpx$ and $\phi_B:B\to\Cpx$ be faithful states.
By adding a scalar multiple of the identity, if necessary, we may without loss of generality assume $\phi_B(b)=0$.
Let
\[
(A*^\red_\Cpx B,\phi)=(A,\phi_A)*_\Cpx(B,\phi_B)
\]
be the reduced free product of C$^*$--algebras.
Using arguments and notation as in the proof of Proposition~\ref{prop:ADBnoninj}, the closure of the
subspace of $L^2(A*_\Cpx^\red B,\phi)$ spanned by elements of the form $(aba')\hat{\;}$ for $a,a'\in A$ is isomorphic to
$\HEu_A\otimes(\Cpx\bh)\otimes\HEu_A$.
To show~\eqref{eq:adbne}, it will suffice to show
\[
(a_1\dt)\hat{\;}\otimes\bh\otimes\ah_2\ne\ah_1\otimes\bh\otimes(\dt a_2)\hat{\;}
\]
in $\HEu_A\otimes(\Cpx\bh)\otimes\HEu_A$.
However, this follows from the assumptions.
\end{proof}

From the above proposition, we get the following example, which requires only ``bad'' relations between $A$ and $\Dt$,
not between $B$ and $\Dt$.
\begin{example}\label{ex:noninj2}\rm
Let $D$, $\Dt$, $A$ and $\At$ be as in Example~\ref{ex:noninjective}.
Let $B$ be any unital C$^*$--algebra of dimension greater than $1$ and let $\Bt=B\otimes\Dt$, (for the unique C$^*$--tensor norm).
Then the $*$--homomorphism $\lambda:A*_DB\to\At*_\Dt\Bt$ arising from the inclusions~\eqref{eq:incl3} is not injective.
\end{example}

\begin{remark}\label{rmk:alg}\rm
The problem with injectivity of $\lambda$ in Examples~\ref{ex:noninjective} and~\ref{ex:noninj2} arises already at the
algebraic level 
\begin{equation}\label{eq:algmap}
A*^\alg_DB\to\At*^\alg_\Dt\Bt.
\end{equation}
On the other hand, in Examples~\ref{exs:triv}, we can arrange that the map between algebras~\eqref{eq:algmap} is injective,
while $\lambda$ fails to be injective, e.g.\ by taking $E$ to be a reduced free product.
However, we do not know of an example where $\lambda$ fails to be injective and where the algebraic map~\eqref{eq:algmap} is injective,
but where $A\cap\Dt=D=B\cap\Dt$.
\end{remark}

\section{An application to residual finite dimensionality}\label{sec:rfd}

A C$^*$--algebra is said to be residually finite dimensional
(r.f.d.) if it has a separating family of finite dimensional
$*$--representations. The first result linking full free products
and residual finite dimensionality was M.-D.~Choi's
proof~\cite{Choi80} that the full group C$^*$--algebras of
nonabelian free groups are r.f.d. In~\cite{EL}, Exel and Loring
proved that the full free product of any two r.f.d.\
C$^*$--algebras $A$ and $B$ with amalgamation over either the zero
C$^*$--algebra or over the scalar multiples of the identity (if
$A$ and $B$ are unital) is r.f.d. In~\cite{BD}, N.~Brown and
Dykema proved that a full amalgamated free product of matrix
algebras $M_k(\Cpx)*_DM_\ell(\Cpx)$ over a unital subalgebra $D$
is r.f.d.\ provided that the normalized traces on $M_k(\Cpx)$ and
$M_\ell(\Cpx)$ restrict to the same trace on $D$. In this section,
we observe that by applying Proposition~\ref{prop:fpembed}, one
obtains (as a corollary of the result from~\cite{BD}) the
analogous result for full amalgamated free products of finite
dimensional algebras.

\begin{lemma}\label{lem:rat}
Let $S=\{ x\in \Reals^n \mid Ax=0\}$, where $A$ is an $m\times n$ matrix having only
rational entries. Then vectors having only rational entries are dense in $S$.
\end{lemma}
\begin{proof}
By considering the reduced row--echelon form of A, we see that
there is a basis for $S$ consisting of rational vectors.
\end{proof}

\begin{thm}\label{thm:rfd}
Consider unital inclusions of C$^*$--algebras $A\supseteq
D\subseteq B$ with $A$ and $B$ finite dimensional. Let $A*_DB$ be
the corresponding full amalgamated free product. Then $A*_DB$ is
residually finite dimensional if and only if there are faithful
tracial states $\tau_A$ on $A$ and $\tau_B$ on $B$ whose
restrictions to $D$ agree.
\end{thm}
\begin{proof}
Since every separable r.f.d.\ C$^*$--algebra has a faithful tracial state, the necessity of the existence of $\tau_A$ and $\tau_B$
is clear.

Let us recall some well known facts about a unital inclusion $D\subseteq A$ of finite dimensional C$^*$--algebras
(see e.g.\ Chapter 2 of~\cite{GHJ}).
Let $p_1,\ldots,p_m$ be the minimal central projections of $A$ and $q_1,\ldots,q_n$ the minimal central projections of $D$.
Then the inclusion matrix $\Lambda^A_D$ is a $m\times n$ integer matrix whose $(i,j)$th entry is
$\rank(q_jp_iAq_j)/\rank(q_jD)$, where the rank of a matrix algebra $M_k(\Cpx)$ is $k$.
To a trace $\tau$ on $A$, we associate the column vector $s$ of length $m$ whose $i$th entry is the trace of a minimal
projection in $p_iA$.
Then the restriction of $\tau$ to $D$ has associated column vector $(\Lambda^A_D)^ts$, where the superscript $t$ indicates transpose.

Thus, given $A\supseteq D\subseteq B$ as in the statement of the
theorem, the existence of faithful tracial states $\tau_A$ and
$\tau_B$ agreeing on $D$ is equivalent to the existence of column
vectors $s_A$ and $s_B$, none of whose components are zero, such
that $(\Lambda^A_D)^ts_A=(\Lambda^B_D)^ts_B$, i.e.\
\begin{equation}\label{eq:Lams}
\left[\begin{array}{rr}(\Lambda^A_D)^t,&-(\Lambda^B_D)^t\end{array}\right]\left[\begin{array}{rr}s_A \\ s_B\end{array}\right]=0.
\end{equation}
Supposing now that such traces $\tau_A$ and $\tau_B$ exist,
by Lemma~\ref{lem:rat} there is a solution $\left[\begin{smallmatrix}s_A\\s_B\end{smallmatrix}\right]$
to~\eqref{eq:Lams} whose entries are all strictly positive and rational.
Therefore, the traces $\tau_A$ and $\tau_B$ agreeing on $D$ can be chosen to take only rational values
on minimal projections of $A$ and, respectively, $B$.
Hence there are unital inclusions into matrix algebras,
\[
M_k(\Cpx)\supseteq A\supseteq D\subseteq B\subseteq M_\ell(\Cpx),
\]
so that $\tau_A$ is the restriction of the tracial state on $M_k(\Cpx)$ to $A$ and $\tau_B$ is the restriction
of the tracial state on $M_\ell(\Cpx)$ to $B$.
By Proposition~\ref{prop:fpembed}, $A*_DB$ is a subalgebra of $M_k(\Cpx)*_D M_\ell(\Cpx)$.
By Theorem~2.3 of~\cite{BD}, $M_k(\Cpx)*_D M_\ell(\Cpx)$ is r.f.d.
Therefore, $A*_DB$ is r.f.d.
\end{proof}

\bibliographystyle{plain}

\end{document}